\documentclass[12pt, a4paper]{article}
\usepackage{a4wide,amssymb,amsmath,amscd,}

\newcommand{\Z}{{\mathbb Z}}
\newcommand{\R}{{\mathbb R}}
\newcommand{\C}{{\mathbb C}}

\begin{document}

\topmargin -2pt

\headheight 0pt

\topskip 0mm \addtolength{\baselineskip}{0.20\baselineskip}
\begin{flushright}
{\tt math.QA/0402401} \\
{\tt KIAS-P04007}
\end{flushright}

\vspace{5mm}

\begin{center}
{\large \bf  Theta Vectors and Quantum Theta Functions }\\

\vspace{5mm}

{\sc Ee Chang-Young}\footnote{cylee@sejong.ac.kr}\\
{\it Department of Physics, Sejong University, Seoul 143-747, Korea}\\

\vspace{2mm}

and \\

\vspace{2mm}

{\sc Hoil Kim}\footnote{hikim@knu.ac.kr}\\

{\it Topology and Geometry Research Center, Kyungpook National University,\\
Taegu 702-701, Korea}\\

\vspace{10mm}

\end{center}

\begin{center}
{\bf ABSTRACT}
\end{center}
In this paper, we clarify the relation between Manin's quantum
theta function and Schwarz's theta vector in comparison with the
$kq$ representation, which is equivalent to the classical theta
function, and the corresponding coordinate space wavefunction.
We
first explain the equivalence relation between the classical theta
function and the $kq$ representation in which the translation
operators of the phase space are commuting.
When the translation operators of the phase space are not
commuting, then the $kq$ representation is no more meaningful.
We explain why Manin's quantum theta function obtained via algebra
(quantum tori) valued inner product of the theta vector is a
natural choice for quantum version of the classical theta function
($kq$ representation). We then show that this approach holds for a
more general theta vector with constant obtained from a
holomorphic connection of constant curvature than the simple
Gaussian one used in the Manin's construction.
We further discuss the properties of the theta vector and of the
quantum theta function, both of which have similar symmetry
properties under translation. \\

\vfill


\thispagestyle{empty}

\newpage
\section*{I. Introduction}

Classical theta functions can be regarded as state functions on
classical tori, and have played an important role in the string
loop calculation \cite{jp,gsw}. Its quantum version on the
noncommutative tori has been discussed mainly by Manin
\cite{manin1,manin2,manin3} and Schwarz \cite{schwarz01,ds02}. In
the physics literature it has been discussed in the context of
noncommutative soliton \cite{mm01}.

In noncommutative field theory, one can find nontrivial soliton
solutions in terms of projection operators
\cite{gms00,mm01,ghs01}.
Before this development, Boca \cite{boca99} has constructed
projection operators on the ${\Z}_4$-orbifold of noncommutative
two torus. There it was also shown that these projection operators
can be expressed in terms of the classical theta functions, of
which certain classical commuting variables are replaced with
quantum operators. Hinted from and generalizing the Boca's result,
Manin \cite{manin2,manin3} explicitly constructed a quantum theta
function, the concept of which he introduced previously
\cite{manin1}.
In both Boca's and Manin's constuctions, the main pillars were the
algebra valued inner product that Rieffel \cite{rief88} used in
his classic work on projective modules over noncommutative tori.
One major difference is that in Manin's construction of quantum
theta function, the so-called theta vector that Schwarz introduced
earlier \cite{schwarz01,ds02} was used for the inner product,
while in Boca's construction the eigenfunctions of Fourier
transform were used.


Both the classical theta function \cite{mumford} and the $kq$
representation in the physics literature \cite{zak72,bgz75} have
been known for a long time. The $kq$ representation is a
transformation of a wavefunction on (real $n$-dimensional)
coordinate space to a function on (real $2n$-dimensional) phase
space consisting of (quasi-)coordinates and (quasi-)momenta.
However, the translation operators in the $kq$ representation
acting on the lattice of the phase space are commuting.
When the lattice of the phase space is periodic, one can identify
functions possessing translational symmetry on the lattice with
the classical theta functions on tori.
When the translation operators of the coordinate and momentum
directions are not commuting, the $kq$ representation and the
classical theta function lose their meaning. One has to find other
ways of representing periodic functions on the lattice of the
non-commuting phase space.
When the algebras are noncommutative, algebra valued inner product
is a good fit for constructing operators out of state functions.
In the case at hand, the coordinates of the phase space are
non-commuting and so is the algebra based on them. And the
functions on the non-commuting phase space can be regarded as
operators.

Classical phase space variables are commuting variables, and thus
they can be simply multiplied in front of a state function
(wavefunction). Namely, we can simply put the values of
observables in front of a statefunction.
 However, in the
quantum case, we have to be very careful with observables.
Quantum
observables behave as operators acting on a state and in general
they change the state.
In fact, the theta vector corresponds to a state on a quantum
torus and the quantum theta function defined by Manin
\cite{manin2,manin3} is an operator acting on the states (module)
on a quantum torus.
In quantum mechanics, one can build operators out of state
vectors. In mathematics, this can be carried out via operator
(algebra) valued inner product. Therefore, it is very natural to
use algebra valued inner product to build the quantum theta
functions from the theta vectors over noncommutative tori.
The classical theta function possesses a certain symmetry property
under the lattice translation, and Manin's quantum theta function
is constructed in such a way that this symmetry property is
maintained as a functional relation which the quantum theta
function should satisfy.


In this paper, we first review the classical theta function and
the $kq$ representation briefly and discuss their relationship. We
then proceed to the quantum case and explain why the Manin's
approach based on algebra valued inner product is a natural choice
for quantum extension. As a support for this viewpoint, we show
that the Manin's construction also holds for a more general theta
vector satisfying the holomorphicity condition. Namely, the
quantum theta function built with our new theta vector also
satisfies the Manin's consistency requirement for the
translational symmetry on the quantum lattice.

%
We also discuss how the theta vectors can be regarded as invariant
state vectors under parallel transport over noncommutative tori
equipped with complex structures, while quantum theta functions
can be regarded as observables having translational symmetry on
the quantum lattice.


The organization of the paper is as follows.
 In section II, we review the classical theta function briefly, then
 explain the relationship between the classical theta functions
 and the $kq$ representation.
 In section III, we first review the theta vectors on quantum
 tori, then explain how the concept of Manin's
 quantum theta function emerges from algebra valued inner product of a state function.
 In section IV, we first review Manin's construction of quantum theta function in detail.
 Then, in order to provide a further support for the Manin's approach
we apply it to the case of a more general theta vector with
constant satisfying the holomorphicity conditon, and show that new
quantum theta function also satisfies the Manin's functional
relation for consistency requirement. In section V, we conclude
with discussion.
\\

\section*{II. Classical complex tori and $kq$ representation}\label{kq-ct}

In this section, we discuss the relationship between the classical
theta function and the so-called $kq$ representation
\cite{zak72,bgz75}. We first look into how the classical theta
function emerges from Gaussian function via Fourier-like
transformation. We then show that the transformed function is
exactly equivalent to the $kq$ representation known in the physics
literature.

We now recall the property of classical theta function briefly,
then show how Gaussian function can be transformed into the
classical theta function.
The classical theta function $\Theta$ is a complex valued function
on ${\C}^n$ satisfying the following relation.
\begin{align}
\Theta(z+\lambda')& =\Theta(z) ~~~~~ {\rm for} ~~~ z \in {\C}^n ,
~ \lambda' \in \Lambda', \label{ct1} \\
\Theta(z+\lambda) & = c(\lambda) e^{q(\lambda , z )} \Theta (z)
~~~~~ {\rm for} ~~~ \lambda \in \Lambda , \label{ct2}
\end{align}
where $~ \Lambda' \bigoplus \Lambda \subset {\C}^n ~$ is a
discrete sublattice of rank $2n$ split into the sum of two
sublattices of rank $n$, isomorphic to $~{\Z}^n~$, and $~c:~
\Lambda \rightarrow {\C}~ $ is a map and $~q: ~ \Lambda \times
{\C} \rightarrow {\C}~$ is a biadditive pairing linear in $z$.

The function $\Theta(z, T)$ satisfying (\ref{ct1}) and (\ref{ct2})
is  defined as
\begin{align}
\Theta(z, T) = \sum_{k\in {\Z}^n} e^{ \pi i (k^t T k + 2 k^t z)}
\label{ct3}
\end{align}
where $T$ is a symmetric complex valued $n\times n$ matrix whose
imaginary part is positive definite. Let $f_T (x) $ be a Gaussian
function defined as below using the same $T$ as above.
\begin{align}
f_{T} (x) = e^{\pi i x^t T x} ~~~ {\rm for} ~~ x \in {\R}^n .
\end{align}
Then $\widetilde{f}_{T} (\rho, \sigma) $ is defined as
\cite{schwarz01}
\begin{align}
\widetilde{f}_{T}(\rho, \sigma) \equiv \sum_{k\in {\Z}^n} e^{-2
\pi i \rho^t k} f_{T} (\sigma + k) \label{ftilde}
\end{align}
where  $~ \rho, \sigma \in {\R}^n$. When we fix $\sigma$, this is
a Fourier transformation between $k$ and $\rho$.
Then from (\ref{ftilde}), we get $\Theta (z,T)$ with a
substitution $z= T \sigma - \rho $ as follows. 
\begin{align}
\widetilde{f}_{T}(\rho, \sigma) & = \sum_{k\in {\Z}^n}  e^{ \pi i
((\sigma + k)^t T (\sigma + k) - 2 \rho^t k)} \label{fkq}
\\
& = e^ { \pi  i \sigma^t T \sigma}   \sum_{k\in {\Z}^n}  e^{ \pi i
( k^t T k + 2 k^t (T \sigma - \rho ))} \nonumber
\\
& = e^ { \pi i \sigma^t T \sigma}  \Theta (T \sigma - \rho, T)
\end{align}

We can do the same procedure for a general Gaussian function,
$~f_{T,c}(x)$, as follows.
\begin{align}
f_{T,c} (x) = e^{\pi i (x^t T x + 2 c^t x)} \label{fc}
\end{align}
where $c \in {\C}^n$. %
Then,
\begin{align}
\widetilde{f}_{T,c}(\rho, \sigma) & \equiv \sum_{k\in {\Z}^n}
e^{-2
\pi i \rho^t k} f_{T,c} (\sigma + k) \label{fctilde} \\
& = \sum_{k\in {\Z}^n}  e^{ \pi i ((\sigma + k)^t T (\sigma + k)
+2 c^t (\sigma +k) - 2 \rho^t k)} \label{fckq}
\\
& = e^ { \pi i ( \sigma^t T \sigma + 2 c^t \sigma ) } \sum_{k\in
{\Z}^n} e^{ \pi  i ( k^t T k + 2 k^t (T \sigma - \rho +c ))}
\nonumber
\\
& = e^ { \pi   i ( \sigma^t T \sigma + 2 c^t \sigma )} \Theta (T
\sigma - \rho +c , T)  .
\end{align}
In this case we get $\Theta (z,T)$ with a substitution $z= T
\sigma - \rho + c .$

The transformation (\ref{ftilde}) exactly matches the
transformation used in defining the $kq$ representation which
already appeared in the physics literature \cite{zak72,bgz75}.
The $kq$ representation is similar to the coherent states for a
simple harmonic oscillator. The coherent states are the
eigenstates of annihilation operator $\hat{a}$, which is a linear
combination of the position and momentum operators. Thus the
eigenvalues of coherent states can be expressed in terms of
expectation values of both position and momentum of the state.
This is in contrast with a usual wavefunction in which position
and momentum eigenvalues do not appear together.

The $kq$ representation which defines symmetric coordinates $k$
(quasimomentum) and $q$ (quasicoordinate) is a transformation from
a wavefunction in position space into a wavefunction in both $k$
and $q$, which we denote as $C(k,q)$.
%
$C(k,q)$ is defined by \cite{bgz75}
\begin{align}
C(k,q) = ( \frac{a}{2 \pi})^{\frac{1}{2}} \sum_{l \in {\Z}}
e^{ikal} \psi (q-la) ~~~ \label{ckq}
\end{align}
where $a$ is a real number (lattice constant), and the
``coordinates" of the phase space ($k,~q$) run over the intervals
$ - \frac{\pi}{a} < k \leqslant \frac{\pi}{a} $ and $ -
\frac{a}{2} < q \leqslant \frac{a}{2} $.
In this representation, the displacement operators $e^{imbx},  ~~
e^{inap}$ in the $x$ and $p$ directions, where $[x,p]=i $,
$b=\frac{2\pi}{a}$, and $ m,n \in {\Z}$, are mutually commuting
and thus they simply become simple multiplication by the function
$e^{im\frac{2 \pi}{a} q}$ and $ e^{inak}$, respectively
\cite{bgz75}.

Comparing (\ref{ckq}) with (\ref{ftilde}), it is not difficult to
see that $C(k,q)$ corresponds to $\widetilde{f}_T (\rho, \sigma)$
in our previous discussion with a correspondence $(\rho
\leftrightarrow k)$ and $( \sigma \leftrightarrow q)$.
Furthermore, from (\ref{ckq}) it can be easily checked that
\begin{align}
& C(k+ \frac{2 \pi}{a}, q)  = C(k,q) , \label{kqprt1} \\
& C(k, q+a)  = e^{ika} C(k,q). \label{kqprt2}
\end{align}
These exactly match (\ref{ct1}) and (\ref{ct2}), the property of
the classical theta function.
We can thus say that the classical theta function corresponds to
the $kq$ representation, $C(k,q)$, while the pre-transformed
Gaussian function $f_T (x)$ for the classical theta function
corresponds to the wavefunction $\psi(x)$ for the $kq$
representation. This correspondence is only valid when the
translation operators of the phase space ($x,~p$) are mutually
commuting.

Therefore, we can see from the above observation that the quantum
theta functions on noncommutative tori cannot be obtained via this
kind of Fourier-like transformation. Since the translation
operators on noncommutative (quantum) tori are in general
non-commuting, we need other ways of going from the position space
representation (like a wavefunction) to the phase space
representation (like $C(k,q)$ or the classical theta function in
the above correspondence) in the quantum case. Namely we have to
find a way to transform a wavefunction (state vector) into an
observable in a noncommuting phase space (consisting of operators
$x$ and $p$). This process can be done via the so-called algebra
valued inner product demonstrated well in the Rieffel's seminal
work on noncommutative tori \cite{rief88}. Manin
\cite{manin2,manin3} has demonstrated sucessfully how this
machinary can be used to define the quantum theta function. We now
turn to this subject in the next section.
\\

%
\section*{III. Theta vectors on quantum tori and algebra valued inner product
           for a passage to quantum theta functions }\label{tv-avip}

In this section, we first discuss theta vectors on quantum tori
and define algebra (quantum tori) valued inner product on the
modules over the quantum tori. Then we introduce Manin's quantum
theta function \cite{manin3} via algebra valued inner product.

A noncommutative $d$-torus $~ T_\theta^d~$
 is a $C^*$-algebra generated by $d$ unitaries $U_1, \dots
, U_d$ subject to the relations
\begin{align}
U_\alpha U_\beta = e^{2 \pi i \theta_{\alpha \beta} } U_{\beta}
U_{\alpha}, ~~~ {\rm for} ~~~ 1\leq \alpha, \beta \leq d,
\end{align}
where $\theta =(\theta_{\alpha \beta} )$ is a skew symmetric
matrix with real entries.

Let $L$ be all derivations on  $T_\theta^d~$, i.e., \[ L = \{
\delta | \delta : ~ T_\theta^d  \rightarrow  T_\theta^d ,~ {\rm
which ~ is ~ linear, ~ and} ~ \delta(fg) =\delta(f)g +f \delta(g)
\} .
\]
Then $L$ has a Lie algebra structure since $[\delta_1 , \delta_2 ]
= \delta_1 \delta_2 - \delta_2 \delta_1  \in L $. We can also see
that $L$ is isomorphic to $ {\R}^d $.
A noncommutative torus is said to have a complex structure if the
Lie algebra $L={\R}^d$ acting on $T_\theta^d$ is equipped with the
complex structure that we explain below. A complex structure on
$L$ can be considered as a decomposition of complexification $L
\bigoplus i L$ of $L$ into a direct sum of two complex conjugate
subspace $L^{1,0}$ and $L^{0,1}$. We denote a basis in $L$ by
$\delta_1, \dots , \delta_d ,$ and a basis in $L^{0,1}$ by
$\tilde{\delta}_1, \dots , \tilde{\delta}_n $ where $ d=2n $.
 One can express
$\tilde{\delta}_\alpha$ in terms of $\delta_j$ as
$\tilde{\delta}_\alpha = t_{\alpha j} \delta_j$, where $t_{\alpha
j}$ is a complex $n\times d$ matrix.

Let $\nabla_j$ (for $j =1, \dots , d$) be a constant curvature
connection on a $T_\theta^d$-module $\cal{E}$. A complex structure
on $\cal{E}$ can be defined as  a collection of ${\C}$ linear
operators $\widetilde{\bigtriangledown}_1, \dots ,
\widetilde{\bigtriangledown}_n $ satisfying
\begin{align}
\widetilde{\bigtriangledown}_\alpha (a \cdot f) & = a
\widetilde{\bigtriangledown}_\alpha f + ( \tilde{\delta}_\alpha a
)
\cdot f \label{holconn} \\
[ \widetilde{\bigtriangledown}_\alpha ,
\widetilde{\bigtriangledown}_\beta ] & =0 \label{comdel}
\end{align}
where $a \in T_\theta^d$ and $f \in \cal{E}$ \cite{schwarz01}. \\
These two conditions are satisfied if we choose
$\widetilde{\bigtriangledown}_\alpha $ as
\begin{align}
\widetilde{\nabla}_\alpha =t_{\alpha j} \nabla_j ~~ {\rm for} ~~
\alpha = 1, \dots, n, ~~j=1, \dots, n .
\end{align}
A vector $f \in \cal{E}$ is holomorphic if
\begin{equation}
\widetilde{\bigtriangledown}_\alpha f =0, ~~~{\rm for}~~~
\alpha=1, \dots , n . \label{holc}
\end{equation}

A finitely generated projective module over $T_\theta^d$ can take
the form $S({\R}^p \times {\Z}^q \times F)$ where $2p + q =d$ and
$F$ is a finite Abelian group \cite{rief88}. Here, $S(M)$ denotes
the Schwartz functions on $M$ which rapidly decay at infinity.

Here, we consider the case that the module is given by
$S({\R}^n)$,
and choose a constant curvature connection $\nabla$
on $S({\R}^n)$ such that
\begin{equation}
( \bigtriangledown_\alpha, \bigtriangledown_{n +\alpha} )  = (
\frac{\partial}{\partial x^\alpha }, -2 \pi i \sigma_\alpha
x_\alpha ) ~~~{\rm for}~~~ \alpha=1, \dots , n , \label{conn}
\end{equation}
where $\sigma_\alpha$ are some real constants, $x^\alpha$ are
coordinate functions on $\mathbb{R}^n$ and repeated indices are
not summed. Then the curvature $[ \bigtriangledown_i ,
\bigtriangledown_j ]=F_{ij}$ satisfies $F_{\alpha, n+\alpha}=2 \pi
i \sigma_\alpha,~F_{n+\alpha, \alpha}= - 2 \pi i \sigma_\alpha$
and all others are zero.
Now, we change the coordinates such that $t =( t_{\alpha j})$
becomes
\begin{align}
t=( {\bf 1}, \tau ),
\end{align}
where ${\bf 1}$ is an identity matrix of size n and $\tau$ is an
$n \times n$ complex valued matrix.

Then, the holomorphic vector $f$ satisfying (\ref{holc}) can be
expressed as
\begin{align}
(\frac{\partial}{\partial x^\alpha} - \sum_{\beta} 2  \pi i
T_{\alpha \beta} x^\beta ) f = 0 ,
\end{align}
where the $n\times n$ matrix $T= (T_{\alpha\beta})$ is given as
follows. The condition (\ref{comdel}) requires that the matrix $T$
be symmetric, $T_{\alpha \beta} = T_{\beta \alpha}$, and it is
given by $T_{\alpha \beta}= \tau_{\alpha \beta} \sigma_\beta, ~
\alpha, \beta = 1, \dots, n $, with the repeated index $\beta$ not
summed.
Up to a constant we get,
\begin{align}
f(x^1, \dots , x^n) =  e^{ \pi i x^\alpha T_{\alpha \beta} x^\beta
}.
\end{align}
If ~${\rm Im} T$ is positive definite, then $f$ belongs to
$S(\mathbb{R}^n)$.
The vectors satisfying the holomorphicity condition (\ref{holc})
are called the theta vectors \cite{schwarz01}.

If a constant in ${\C}^n $ is added to a given connection
$\widetilde{\bigtriangledown}$, it still yields the same constant
curvature. Then the holomorphicity condition (\ref{holc}) becomes
\cite{ds02,kl03}
\begin{equation}
(\widetilde{\bigtriangledown}_\alpha -2 \pi i c_\alpha ) f_c =0,
~~~{\rm for}~~~ \alpha=1, \dots , n  \label{holcc}
\end{equation}
 for $f_c \in S({\R}^n)$, giving the following condition
\begin{align}
(\frac{\partial}{\partial x^\alpha} -\sum_{\beta}  2  \pi i
T_{\alpha \beta} x^\beta - 2  \pi i c_\alpha) f_c = 0 ,
\end{align}
whose solution we get
\begin{align}
 f_c(x) = e^{ \pi i x^\alpha T_{\alpha \beta} x^\beta + 2 \pi i c_\alpha x^\alpha
 }.
\end{align}
Here, we would like to make an observation.
 The holomorphicity condition (\ref{holc}) means that
the theta vector $f $ or $f_c$ is invariant under a parallel
transport on a noncommutative torus with complex structure.

Now we turn to the concept of the quantum theta function
introduced by Manin \cite{manin1,manin2,manin3}. Recall that the
classical theta function $\Theta(z)$ satisfies the conditions
(\ref{ct1}) and (\ref{ct2})
\begin{align*}
 \Theta(z+\lambda') & = \Theta(z), ~~~ z \in {\C}^n, ~~ {}^\forall  \lambda' \in \Lambda'
 ,  \\
 \Theta(z+\lambda) & = c(\lambda) e^{q(\lambda,z)} \Theta(z), ~~~ {}^\forall \lambda \in \Lambda
 ,
\end{align*}
where $~ c: ~ \Lambda \rightarrow {\C} ~$ is a map and $~ q:~
\Lambda\times {\C} \rightarrow {\C}~$ is a biadditive pairing
linear in $z$. This function can be written formally as follows
\cite{manin1}.
\begin{align}
\Theta(z) = \sum_{j \in J} a_{j} e^{ 2 \pi i j(z)},
\end{align}
where $J= {\rm Hom}(\Lambda' , {\Z})$. The coefficients $a_j$
decay swiftly enough. Then this form satisfies the first condition
(\ref{ct1}) automatically and we impose a constraint for $a_j$
satisfying the second condition (\ref{ct2}). If we define
$T(J)({\C}) = {\rm Hom} (J, {\C}^*)$ where ${\C}^* = {\C} - \{ 0
\}$.
 We have an
isomorphism $e$ from $J$ to $ \widetilde{J} \equiv {\rm Hom}
(T(J)({\C}), {\C}^*)  $.
We denote $e(j)$ the image of $j$ by this map $e$. Then
\begin{equation*}
e(j+l)=e(j) e(l), ~~{\rm for}~~ j,l \in J .
\end{equation*}
We have an analytic map $P$ which is in fact an isomorphism up to
$\Lambda'$,
\begin{equation*}
P : ~ {\C}^n \longrightarrow T(J)({\C}) ,
\end{equation*}
inducing the pullback $ P^*(e(j)) = e^{2 \pi i j(\cdot)}  $ where
$j(\cdot)$ is the linear function on ${\C}^n$ extending $j$ as a
function on $ \Lambda'$.
Then the classical theta function $\Theta$ can be expressed as
\[ \Theta = P^*(\widetilde{\Theta}), ~~{\rm where}~~
    \widetilde{\Theta} = \sum_{j \in J} a_j e(j) . \]
Let $B$ be the image of $\Lambda$ under $P$, then
$b^*(\widetilde{\Theta})$, the translation of $\widetilde{\Theta}$
by $b \in B$, is equal to $\sum_{j \in J} a_j j(b) e(j)$, where
$j(b)=e(j)(b)$ is the value of $e(j)$ at the point $b \in B$:
\[ b^*(\widetilde{\Theta})(w) = \widetilde{\Theta}(w \cdot b),
~~ {\rm where} ~~ \forall w \in T(J)({\C}) . \]
The second condition can be interpreted as
\begin{align}
c_b e(j_b) b^*(\widetilde{\Theta}) = \widetilde{\Theta}
\label{autofac}
\end{align}
where $c_b \in {\C}$ and $j_b  \in J$.
 To generalize this for $T_\theta^d$, the
Heisenberg group $G(J)$ is defined. This is the group of linear
endomorphisms of the space of functions $(\Phi)$ on algebraic
torus $T(J)({\C})$ generated by the following maps,
\begin{align}
 [ c, x, j  ] : ~ \Phi \rightarrow c e(j) x^* (\Phi ) ,
\label{trsmap}
\end{align}
where $c\in {\C}^*, ~ x \in T(J) ({\C}),~ j \in J $ and
$x^*(e(j))=j(x) e(j)$, where $j(x)$ being the value of $e(j)$ at
$x$. In these terms, a system consisting of a subgroup $B$ in
$T(J)({\C})$ and automorphy factors satisfying the second
condition (\ref{autofac}) become simply a homomorphism, which we
will call a multiplier, $\cal L$,
\begin{align}
{\cal L}: ~ B \rightarrow G(J), ~ {\cal L}(b) = [ c_b, x_b, j_b ],
\label{automorphy}
\end{align}
where  $b \rightarrow x_b $ is a bijection.
Manin's quantum theta function is invariant under the image of
${\cal L}$, the subgroup of the Heisenberg group $G(J)$.

Now, we consider the algebra valued inner product on a bimodule
after Rieffel \cite{rief88}.
Let $M$ be any locally compact Abelian group, and $\widehat{M}$ be
its dual group and ${\cal G} \equiv M \times \widehat{M} $.
Let $\pi$ be a representation of ${\cal G}$ on $L^2(M)$ such that
\begin{align}
\pi_x \pi_y = \alpha (x,y) \pi_{x+y} =\alpha (x,y)
\overline{\alpha}(y,x) \pi_y \pi_x ~~~ {\rm for}~~ x,y \in {\cal
G} \label{ccl}
\end{align}
where $\alpha$ is a map $ \alpha : ~ {\cal G} \times {\cal G}
\rightarrow {\C}^* $ satisfying
\[ \alpha(x,y)
=\alpha(y,x)^{-1} , ~~~ \alpha(x_1 + x_2 , y) = \alpha(x_1 , y)
\alpha (x_2 , y) ,  \] and $\overline{\alpha}$ denotes the complex
conjugation of $\alpha$.

 Let $D$ be a
discrete subgroup of ${\cal G}$. We define $S(D)$ as the space of
Schwartz functions on $D$.
 For $\Phi \in S(D)$, it can be expressed as $\Phi = \sum_{w \in D} \Phi(w) e_{D,
 \alpha}(w)$ where $e_{D, \alpha}(w)$ is a delta function with
 support at $w$ and obeys the following relation.
\begin{equation}
e_{D, \alpha} (w_1) e_{D, \alpha} (w_2) = \alpha(w_1,w_2) e_{D,
\alpha} (w_1 +w_2) \label{ccld}
\end{equation}
For Schwartz functions $f,g \in S(M)$, the algebra ($S(D)$) valued
inner product is defined as
\begin{align}
{}_D <f,g> \equiv \sum_{w\in D} {}_D<f,g>(w) ~ e_{D, \alpha}(w) ~
\label{aip}
\end{align}
where
\begin{align}
{}_D<f,g>(w) = <f, \pi_w g> . \nonumber
\end{align}
Here, the scalar product of the type $<f,p>$ used above for $f,p
\in L^2 (M)$ denotes the following.
\begin{align}
<f,p> = \int f(x_1) \overline{p(x_1)} d \mu_{x_1}  ~~~{\rm for} ~~
x=(x_1,x_2) \in M \times \widehat{M} , \label{sp}
\end{align}
where $\mu_{x_1}$ represents the Haar measure on $M$ and
$\overline{p(x_1)}$ denotes the complex conjugation of $p(x_1)$.
 Thus
the $S(D)$-valued inner product can be represented as
\begin{align}
{}_D <f,g> =\sum_{w\in D} <f, \pi_w g> ~ e_{D, \alpha}(w) ~.
\label{aipr}
\end{align}
%

 For $\Phi \in S(D) $ and $f \in
S(M)$, then $~ \pi (\Phi) f \in S(M)$ can be written as
\cite{rief88}
\begin{align}
(\pi(\Phi)f)(m) & = \sum_{w \in D} \Phi (w) (\pi_w f) (m)
\end{align}
where $m\in M, ~ w=(w', w'') \in D \subset M \times \widehat{M}$.
 ~ For $ f,g \in S(M)$ and $\Phi \in S(D)$, one can also check the
following relation \cite{rief88}
\begin{align}
{}_D< \Phi f, g > = \Phi * {}_D <f,g> ,
\end{align}
where $*$ denotes the convolution.
This means the compatibility of
the $S(D)$-valued inner product with the action of $S(D)$ on
$S(M)$.
%
Now one can define $D^\bot$, the set of $z$'s in ${\cal G}$ such
that $\pi_z$ commutes with $\pi_w$ for all $w\in D$,
 \[ D^\bot =\{ z \in {\cal G} : \alpha(w,z)
\overline{\alpha} (z,w) = 1 , ~~ {}^\forall w \in D \}.
\]
 Then the action of $\Omega \in S(D^\bot)$ on $f \in S(M)$ can be defined as,
\begin{align}
f \Omega = \sum_{z \in D^\bot} (\pi_z^* f) \Omega(z) ,
\end{align}
and thus the $S(D^\bot)$-valued inner product can be expressed as
\begin{align}
<f,g>_{D^\bot} &= \sum_{z \in {D^\bot}}e_{D, \alpha}^*(z)
<f,g>_{D^\bot}(z)
\nonumber \\
 & = \sum_{z \in {D^\bot}}e_{D, \alpha}^*(z) < \pi_z g, f> ,
\end{align}
where $*$ denotes the adjoint operation.
  From the above definitions, the following relation holds \cite{rief88}.
\begin{align}
{}_D <f,g>h =f <g,h>_{D^\bot} ~~~ {\rm for} ~~~ f,g,h \in S(M) .
\end{align}
  Furthermore, if $<f,f>_{D^\bot} = 1$, then ${}_D<f,f>$ is a
projection operator \cite{rief88,manin2,manin3}.

The Manin's quantum theta function $\Theta_D$ \cite{manin2,manin3}
was defined via algebra valued inner product up to a constant
factor,
\begin{align}
{}_D<f_{T} , f_{T} > &  \sim  \Theta_D ,
\label{qtheta-def}
\end{align}
where $f_T$ used in the construction was a simple Gaussian theta
vector
\begin{align}
f_T = e^{\pi i x_1^t T x_1}, ~~ x_1 \in M, \label{tv-gauss}
\end{align}
with $T$ be an $n \times n$ complex valued matrix.
Manin required that the quantum theta function $\Theta_D$ defined
in this way should satisfy the following condition under
translation derived from the map (\ref{trsmap})
\begin{equation}
{}^\forall g \in D, ~~ C_g ~ e_{D, \alpha} (g) ~ x_g^* ( \Theta_D)
= \Theta_D \label{qt-trs}
\end{equation}
where $C_g$ is an appropriately given constant, and $x_g^* $ is a
``quantum translation operator" defined as
\begin{align}
x_g^* (e_{D, \alpha} (h)) =  {\cal X}(g,h) e_{D, \alpha} (h)
\label{qtr-ftn}
\end{align}
with some commuting function  ${\cal X}(g,h)$ for $ g,h \in D$.
 The requirement (\ref{qt-trs})
can be regarded as the quantum counterpart of the second property
of the classical theta function, (\ref{ct2}).

In physics language, the theta vector corresponds to a state
vector (wavefunction) which can be expressed as a Dirac ket, say
$|n>$, and the quantum theta function corresponds to an operator
for an observable which in terms of the Dirac bra-ket notation can
be represented as $\sum_n a_n |n><n|$ with $ a_n \in {\C}$. In the
case of algebra valued inner product, ${}_D<f,f>$ corresponds to
$\sum_n a_n |n><n| ~ \ncong  ~ {\bf 1}$, and $<f,f>_{D_{\bot}} $
corresponds to a case in which $ \sum_n a_n <n|n> ~  \cong  ~ {\bf
1}$.
Namely, the inner product in the latter case becomes a scalar
which is equivalent to an identity operator. Furthermore, as we
mentioned above, (\ref{qt-trs}) represents the quantum version of
the symmetry of the classical theta function under translation.
Thus based on our above discussion in the Dirac's notation and the
symmetry property that we mentioned, we can deduce that the
Manin's quantum theta function constructed via algebra valued
inner product is the quantum version of the classical theta
function.
\\

%
%

\section*{IV. Quantum theta functions - extended to holomorphic connections with constants}\label{qtfc}

In this section, we review Manin's construction of quantum theta
function in detail starting from the algebra valued inner product
of the Gaussian theta vector, and show that Manin's approach for
quantum theta function also holds for the case of a theta vector
obtained from more general holomorphic connections with constants.

As in the classical theta function case, we first introduce an
$n$-dimensional complex variable $\underline{x} \in {\C}^n$ with
complex structure $T$ explained in the previous sections as
\begin{align}
\underline{x} \equiv T x_1 +x_2
\end{align}
where $x=(x_1, x_2) \in M \times \widehat{M}$.
 Based on the defining
concept for quantum theta function (\ref{qtheta-def}), Manin
defined the quantum theta function $\Theta_D$ as
\begin{align}
{}_D<f_{T} , f_{T} > & = \frac{1}{\sqrt{2^n \det ({\rm Im} ~ T )}}
\Theta_D \label{qtfM}
\end{align}
with $f_T$ given by (\ref{tv-gauss}).
Using (\ref{aip}) the $S(D)$-valued inner product in
(\ref{qtheta-def}) can be expressed as
\begin{align}
{}_D<f_{T} , f_{T} >  =\sum _{h \in D} <f_{T} , \pi_h f_{T}> e_{D,
\alpha} (h) . \label{sdip}
\end{align}
Now, we define $\pi$ of ${\cal G}$ on $L^2(M)$ as follows.
\begin{align}
(\pi_{(y_1, y_2)} f)(x_1) =  e^{2 \pi i x_1^t y_2 +\pi i y_1^t y_2
} f(x_1 + y_1 ), ~~~ {\rm for} ~~  x,y \in {\cal G} =M \times
\widehat{M}  \label{uaction}
\end{align}
Then the cocycle $\alpha(x,y)$ in (\ref{ccl}) is given by $
\alpha(x, y) = e^{\pi i (x_1^ty_2 -y_1^t x_2)} $.
%
%

  In \cite{manin3}, Manin  showed that the quantum theta
function defined in (\ref{qtfM}) is given by
\begin{align}
 \Theta_D & = \sum _{h \in D}  e^{- \frac{\pi}{2} H(\underline{h},\underline{h}) }
  e_{D, \alpha} (h) ,
\label{TDM}
\end{align}
where
\[
H( \underline{g}, \underline{h} ) \equiv \underline{g}^t ( {\rm
Im} T)^{-1} \underline{h}^*
\]
with $ \underline{h}^* = \overline{T} h_1 + h_2 $ denoting the
complex conjugate of $\underline{h}$, and satisfies the following
functional equation.
\begin{equation}
{}^\forall g \in D, ~~ C_g ~ e_{D, \alpha} (g) ~ x_g^* ( \Theta_D)
= \Theta_D
\label{TDfnr}
\end{equation}
where $C_g$ is defined by
\[ C_g = e^{- \frac{\pi}{2} H(\underline{g},\underline{g})} \]
and the action of ``quantum translation operator"  $x_g^*$ is
given by
\begin{align}
x_g^* (e_{D, \alpha} (h)) = e^{- \pi
H(\underline{g},\underline{h})} e_{D, \alpha} (h).
\label{xtrans}
\end{align}

We now sketch the proof of the above statement.
The scalar product inside the summation in (\ref{sdip}) can be
expressed as
\begin{align}
<f_{T} , \pi_h f_{T} >  = \int_{{\R}^n} d\mu_{x_1} e^{ \pi i x_1^t
T x_1 -\pi i (x_1 +h_1)^t \overline{T} (x_1+h_1) - 2 \pi i x_1^t
h_2  -\pi ih_1^t h_2} .
\end{align}
Denoting the exponent inside the integral sign  as
\[ e^{-\pi (q(x_1)+ l_h (x_1)+ \widetilde{C}_h)}  \]
with
\begin{align*}
q(x_1) & = 2x_1^t ~ ({\rm Im} T) ~x_1 \\
l_h(x_1) & = 2 i x_1^t (\overline{T}h_1 + h_2 ) \\
\widetilde{C}_h & =i h_1^t (\overline{T}h_1 + h_2 ) ,
\end{align*}
and using the relation
 \[  q(x_1 + \lambda_h) -q(\lambda_h) = q(x_1) + l_h (x_1)  \]
with
 \[ \lambda_h \equiv \frac{i}{2} ( {\rm Im} T)^{-1} \underline{h}^* , \]
the integration now becomes
\begin{align*}
  \int_{{\R}^n}d \mu_{x_1}  e^{- \pi (q(x_1) + l_h(x_1) + \widetilde{C}_h )}  =
  e^{- \pi ( \widetilde{C}_h - q(\lambda_h))}
 \int_{{\R}^n}d \mu_{x_1}  e^{- \pi q(x_1 + \lambda_h)}
 = \frac{1}{\sqrt{\det q}} e^{- \pi (\widetilde{C}_h -q(\lambda_h))} .
\end{align*}
With a straightforward calculation one can check that
\[ \widetilde{C}_h -q(\lambda_h) = \frac{1}{2} H(\underline{h},
\underline{h}) , \]
 and with $ \det q=2^n \det ({\rm Im} ~ T) $,
the expression for Manin's quantum theta function (\ref{TDM})
follows.

The functional relation for quantum theta function (\ref{TDfnr})
 can be shown by use of the definition of ``quantum translation
 operator"  (\ref{xtrans}) as follows. 
\begin{align*}
& C_g ~ e_{D, \alpha} (g) ~ x_g^* (\sum _{h \in D}  e^{-
\frac{\pi}{2} H(\underline{h},\underline{h}) }
  e_{D, \alpha} (h) ) \\
 & = e^{- \frac{\pi}{2} H(\underline{g},\underline{g}) }
  e_{D, \alpha} (g) \sum _{h \in D}  e^{-
\frac{\pi}{2} H(\underline{h},\underline{h})- \pi
H(\underline{g},\underline{h}) }
  e_{D, \alpha} (h)  \\
 & = \sum _{h \in D}  e^{-
\frac{\pi}{2} H(\underline{g} + \underline{h},\underline{g} +
\underline{h})}
  e_{D, \alpha} (g+h)
\end{align*}
In the last step, the cocycle condition (\ref{ccld}) with $
\alpha(g, h) = e^{\pi i (g_1^t h_2 -h_1^t g_2)} = e^{ \pi i {\rm
Im} H(\underline{g},\underline{h}) }$
 was used. This proves the
statement. $\Box$
%

%
%

%

%
%

In the rest of this section, we apply the Manin's approach to a
more general theta vector with constant obtained from a
holomorphic connection of constant curvature. We do this to
provide a further support for Manin's quantum theta function
approach based on the algebra valued inner product and to show
that it is a natural choice for quantum extension of the classical
theta function.

We begin again with $S(D)$-valued inner product (\ref{qtheta-def})
 with a more
general theta vector $f_{T,c} ~$ which appeared in
\cite{ds02,kl03}.
\begin{align}
{}_D<f_{T,c} , f_{T,c} > & =\sum _{h \in D} <f_{T,c} , \pi_h
f_{T,c}> e_{D, \alpha} (h) \label{aip-c}
\end{align}
where
\begin{equation}
f_{T,c} (x_1)=e^{\pi i x_{1}^t T x_1 + 2 \pi i c^t x_{1}}, ~~ c
\in {\C}^n , ~~ x_1 \in M ,
\end{equation}
and $T$ is the complex structure mentioned before.
 From (\ref{sp}) and (\ref{uaction}), the algebra valued inner
product (\ref{aip-c}) can be written as
\begin{align}
{}_D<f_{T,c} , f_{T,c} > & =\sum _{h \in D} <f_{T,c} ,
\pi_h f_{T,c} >  e_{D, \alpha} (h) \nonumber \\
 & =\sum _{h \in D}
 \int_{{\R}^n}d \mu_{x_1}
f_{T,c}(x_1)
\overline{( \pi_h f_{T,c} )(x_1)}  e_{D, \alpha} (h) \nonumber \\
& \equiv \sum _{h \in D} \int_{{\R}^n}d \mu_{x_1} e^{- \pi [q(x_1)
+ l_{h,c} (x_1) +
 \widetilde{C}_{h,c} ] }
 e_{D, \alpha} (h)
\label{aip-c2}
\end{align}
where $q(x_1),~ l_{h,c} (x_1) , ~
 \widetilde{C}_{h,c} $ are defined by
\begin{align}
q(x_1) = & 2 x_1^t ({\rm Im}~ T) x_1 , \nonumber \\
 l_{h,c} (x_1) = & 2 i x_1^t (\overline{T} h_1 + h_2 - 2 i ({\rm Im}~ c)), \\
 \widetilde{C}_{h,c} = & i h_1^t (\overline{T} h_1 + h_2 + 2 \overline{c}) . \nonumber
\end{align}
%
%
Denoting
 \[ \lambda_{h,c} \equiv \frac{i}{2} ({\rm Im}~ T)^{-1} (\underline{h}^* - 2 i ({\rm Im}~ c)), \]
 one can check that
\[ q(x_1) + l_{h,c} (x_1)  = q(x_1 + \lambda_{h,c}) -q(\lambda_{h,c}) . \]
%
Thus, the algebra valued inner product (\ref{aip-c2}) can be
written as
\begin{align}
{}_D<f_{T,c} , f_{T,c} > & =\sum _{h \in D}  e^{- \pi (
\widetilde{C}_{h,c} - q(\lambda_{h,c}))} e_{D, \alpha} (h)
 \int_{{\R}^n}d \mu_{x_1}  e^{- \pi q(x_1 + \lambda_{h,c})} .
\end{align}
Since $\int_{{\R}^n}d \mu_{x_1}  e^{- \pi q(x_1 + \lambda_{h,c})}=
1/ \sqrt{\det q} $, the above expression can be rewritten as
\begin{align}
{}_D<f_{T,c} , f_{T,c} > & = \frac{1}{\sqrt{2^n \det ({\rm Im} ~ T
)}} \sum _{h \in D}  e^{- \pi ( \widetilde{C}_{h,c} -
q(\lambda_{h,c}))} e_{D, \alpha} (h)
\end{align}
and we define our quantum theta function $\Theta_{D,c}$ as
\begin{align}
{}_D<f_{T,c} ~,~ f_{T,c} > & \equiv \frac{1}{\sqrt{2^n \det ({\rm
Im} ~ T )}} \Theta_{D,c} ~ .
\end{align}
The quantum theta function defined above is evaluated as
\begin{align}
 \Theta_{D,c} & = \sum _{h \in D}  e^{- \pi ( \widetilde{C}_{h,c} -
q(\lambda_{h,c}))} e_{D, \alpha} (h) \nonumber \\
& = \sum _{h \in D}  e^{- \pi [ \frac{1}{2} (\underline{h}^t - 2 i
({\rm Im}~ c)^t) ({\rm Im}~ T)^{-1} (\underline{h}^* - 2 i ({\rm
Im}~ c)) +2 i h_1^t ({\rm Re}~ c)] } e_{D, \alpha} (h).
\label{qthc}
\end{align}
And the above defined quantum theta function $\Theta_{D,c}$
satisfies the following.
\\
%
{\bf Theorem:}
{\it The quantum theta function $\Theta_{D,c}$
defined by the following algebra valued inner product
\begin{align}
{}_D<f_{T,c} ~,~ f_{T,c} > & \equiv \frac{1}{\sqrt{2^n \det ({\rm
Im} ~ T )}} \Theta_{D,c} ~
\end{align}
with a theta vector $f_{T,c}$  below, which is obtained from a
holomorphic connection with constant $c \in {\C}^n$,
\begin{equation}
f_{T,c} (x_1)=e^{\pi i x_{1}^t T x_1 + 2 \pi i c^t x_{1}},
\end{equation}
satisfies the following identity
\begin{equation}
{}^\forall g \in D, ~~ C_{g,c} ~ e_{D, \alpha} (g) ~ x_{g,c}^* (
\Theta_{D,c}) = \Theta_{D,c } ~ . \label{qthcfn}
\end{equation}
Here $C_{g,c}  $ is a constant defined by
\[ C_{g,c} \equiv e^{- \frac{\pi}{2} H_c(\underline{g},\underline{g})} \]
where $H_c(\underline{g},\underline{g})$ is given by
\begin{align}
 H_c(\underline{g},\underline{g})
   & = (\underline{g} -2i   ( {\rm Im}~ c  ) )^t
         ({\rm Im}~ T)^{-1} (\underline{g}^* -2i   ( {\rm Im}~ c  ) )
          +  4i g_1^t ( {\rm Re}~ c) ,
          \label{hc}
\end{align}
and  $ ~ x_{g,c}^* ~$ is a ``quantum translation operator" defined
by
\begin{align}
x_{g,c}^* (e_{D, \alpha} (h)) \equiv e^{- \pi
X(\underline{g},\underline{h})} e_{D, \alpha} (h)
\end{align}
where $X(\underline{g},\underline{h})$ is given by }
\[ X(\underline{g},\underline{h}) =  \underline{g}^t ({\rm Im}~ T)^{-1}
\underline{h}^* + 2 ({\rm Im}~ c)^t ({\rm Im}~ T)^{-1} ({\rm Im}~
c).
\]
%
%
{\bf Proof.}
We first note that from (\ref{qthc}) and (\ref{hc}) our quantum
theta function $\Theta_{D,c}$ can be expressed as
\begin{align}
 \Theta_{D,c} & = \sum _{h \in D}  e^{- \frac{\pi}{2}
 H_c(\underline{h},\underline{h})}  e_{D, \alpha} (h) .
\label{qthch}
\end{align}
Thus the left hand side of the functional relation (\ref{qthcfn})
can be written as
\begin{align*}
~ C_{g,c} ~ e_{D, \alpha} (g) ~ x_{g,c}^* ( \Theta_{D,c})
 & = e^{- \frac{\pi}{2}
  H_c(\underline{g},\underline{g})}
e_{D, \alpha} (g) ~ x_{g,c}^*
 (\sum _{h \in D}
  e^{- \frac{\pi}{2} H_c(\underline{h},\underline{h})}
e_{D, \alpha} (h))
\\
  & = \sum _{h \in D} e^{- \frac{\pi}{2}
  H_c(\underline{g},\underline{g})}
e^{- \frac{\pi}{2} H_c(\underline{h},\underline{h})} e_{D, \alpha}
(g) ~ x_{g,c}^* (e_{D, \alpha} (h))
 \\
& = \sum _{h \in D} e^{- \frac{\pi}{2}
  H_c(\underline{g},\underline{g})}
e^{- \frac{\pi}{2} H_c(\underline{h},\underline{h})} e^{- \pi
X(\underline{g},\underline{h})} e_{D, \alpha} (g)
 e_{D, \alpha} (h).
\end{align*}
Then  using the cocycle relation (\ref{ccld})
\[ e_{D, \alpha} (g) e_{D, \alpha} (h) = \alpha(g,h) e_{D, \alpha}
(g+h) = e^{ \pi i {\rm Im}( \underline{g}^t ({\rm Im}~ T)^{-1}
\underline{h}^* ) } e_{D, \alpha} (g+h)  , \] one can check that
with a straightforward calculation
\[
e^{- \frac{\pi}{2}
  H_c(\underline{g},\underline{g})}
e^{- \frac{\pi}{2} H_c(\underline{h},\underline{h})} e^{- \pi
X(\underline{g},\underline{h})} e^{ \pi i {\rm Im}(
\underline{g}^t ({\rm Im}~ T)^{-1} \underline{h}^* ) } = e^{-
\frac{\pi}{2}  H_c(\underline{g}+ \underline{h} ,\underline{g}+
\underline{h})} ,
\]
proving the relation (\ref{qthcfn}).  $ \Box $
\\

The property of quantum theta function (\ref{qthcfn}) represents
the translational symmetry of the quantum lattice. This
corresponds to the symmetry property (\ref{ct2}) of the classical
theta function on the complex tori:
\begin{equation*}
\Theta(z+\lambda) = C(\lambda) e^{q(\lambda,z)} \Theta(z) ~~~{\rm
for }~~~ \lambda \in \Lambda
\end{equation*}
where $\Lambda$ is the period lattice for the complex tori.
The relation is the same as in the case of Manin's construction
expressed in (\ref{TDfnr}). The only difference here is that the
constant factor $C_g$ and the action of ``quantum translation
operator" $x_g^*$ have been changed slightly due to the constant
$c \in {\C}^n$ appearing in our new theta vector $f_{T,c}$. The
changes in these two were possible due to quantum nature of the
quantum theta functions which inherit the mapping property
(\ref{trsmap}) expressed as  a multiplier $\cal{L}$ in
(\ref{automorphy}). For the multiplier $\cal{L}$,
  we have a freedom to
select $c_b$ and $j_b$ in (\ref{automorphy}). The constant factor
$C_g$ and the action of ``quantum translation operator" $x_g^*$
directly corresponds and is related to $c_b$ and $j_b$,
 respectively.
\\

\section*{V. Conclusion }

In this paper we explained how Manin's quantum theta functions
emerge naturally from the state vectors on quantum
(noncommutative) tori via algebra valued inner product.


As we discussed in section III, the theta vectors can be regarded
as invariant state vectors under parallel transport on the
noncommutative tori equipped with complex structures. However,
they are not like the classical theta functions which are the
state vectors (holomorphic sections of line bundles) over
classical tori. This is because the classical theta functions
(complex $n$ dimensional) are equivalent to $kq$ representations
(real $2n$ dimensional) which are transformations of the functions
over coordinates (real $n$ dimensional) only. Namely, these are
functions over the phase space (real $2n$ dimensional) consisting
of coordinates and their canonical momenta, while the theta
vectors are more or less corresponding to the functions over
coordinates (real $n$ dimensional) only.

Therefore to build a quantum version of classical theta function,
we need to build a function over the quantum phase space (real
$2n$ dimensional) via a transformation like $kq$ representation.
However, a function over quantum phase space is necessarily an
operator since coordinates and their momenta are not commuting in
general. As we discussed in section III, the algebra valued inner
product is a good fit for this purpose, since it transforms a
(commuting) function into an operator.
Thus the quantum theta function obtained via algebra valued inner
product from the theta vector (a function over commuting
variables) can be regarded as a quantum version of $kq$
representation which corresponds to the classical theta function.

In conclusion, we can say that the quantum theta function is a
quantum version of the classical theta function which is
equivalent to the $kq$ representation, while the theta vector
corresponds to a wavefunction  over commuting coordinates, the
pre-transformed function for the $kq$ representation.

Finally, we compare the characteristics of the quantum theta
function and the theta vector.  The theta vectors can be regarded
as invariant state vectors under parallel transport on the
noncommutative tori equipped with complex structures, since they
are defined to vanish under the action of the holomorphic
connection which can be regarded as the generator for parallel
transport.
While the quantum theta functions can be regarded as observables
having translational symmetry on the quantum lattice. Thus it is
not surprising that these two are related by algebra valued inner
product which one can regard as a quantum version of the
transformation for the $kq$ representation.
\\

%
%


%


\vspace{5mm}

\noindent
{\Large \bf Acknowledgments}

\vspace{5mm}
\noindent
Most part of the work was done during
authors' visit to KIAS. The authors would like to thank KIAS for
its kind hospitality. This work was supported by KOSEF
Interdisciplinary Research Grant No. R01-2000-000-00022-0.


\vspace{5mm}


\end{document}